\documentclass[a4paper, 12pt]{article}
 \usepackage{graphicx}
\usepackage{graphics}
\usepackage{titling}
\usepackage{enumitem}
\usepackage{textcomp}

\usepackage[T1]{fontenc}

\newcommand{\Addresses}{{
  \bigskip
  \footnotesize
\begin{flushright}
Darya E. Apushkinskaya\\
Saarland University\\ P.O. Box 151150,
66041 Saarbr{\"u}cken, Germany;\\
Peoples' Friendship University of Russia (RUDN University)\\
6 Miklukho-Maklaya St, Moscow, 117198, Russia\\
\textit{darya@math.uni-sb.de}\\
\bigskip
Alexander I. Nazarov\\
St. Petersburg Department of Steklov Institute\\ Fontanka 27,
St. Petersburg 191023, Russia;\\
St. Petersburg State University\\
Universitetskii pr. 28,
 St. Petersburg 198504, Russia\\
\textit{al.il.nazarov@gmail.com} \\
\bigskip
Galina I. Sinkevich\\
St. Petersburg State University of Architecture 
and Civil Engineering\\
 2-ya Krasnoarmeyskaya Ulitsa, 4, Sankt-Peterburg, 190005, Russia\\
\textit{galina.sinkevich@gmail.com}
 \end{flushright} 
}}

\begin{document}
\title{In Search of Shadows: \\ the First Topological Conference,\\ Moscow 1935}

\author{D.E.  Apushkinskaya, A.I. Nazarov,  G.I. Sinkevich}

\maketitle

\begin{abstract}
We discuss some mistakes and curiosities concerned with the celebrated First International Topological Conference in Moscow, 1935.
\end{abstract}

\vspace{0.1cm}

The First International Topological Conference took place in Moscow, September 4--10, 1935. In fact, it was the first truly international specialized topological meeting  in the history of the world mathematical community (see \cite{Al80}). The conference brought together many outstanding experts from 10 countries. Moreover, it gave some ``\dots major breakthroughs toward future methods in topology of great import for the future of the subject.'' (\cite{Whi88}). 

It should be emphasized that to organize this Moscow meeting was not an easy task because of the political situation in the USSR. The conference was realized only thanks to the great efforts of P.S.~Aleksandrov. This was reflected, in particular, in Aleksandrov's correspondence with A.N.~Kolmogorov (see \cite{K03}). Shortly after the Moscow meeting the Iron Curtain was dropped completely that separated Soviet science from the world community for more than two decades.

There are several official publications, historical notes and reminiscences devoted to this conference (\cite{Al36}, \cite{Al80}, \cite{Bi89}, \cite{B36}, \cite{Bul35}, \cite{H66}, \cite{J99}, \cite{S17}, \cite{Tuc35}, \cite{Whi88}, and \cite{Z35}). However, it turns out that some data  (number of speakers, number of talks, etc.) diverge and sometimes contradict each other in different sources. 

Based on the materials  at our disposal, we suggest our own version of this story.

In our opinion, the most complete  information  is given in the survey of A.F. Lapko and L.A. Lyusternik \cite[pp.82--85]{LL57}. Namely, the authors of \cite{LL57}  provide the following list of talks:
\vspace{0.2cm}

\begin{small}
\begin{itemize} [nolistsep]
\item[1.] J.W.~Alexander (USA), \textit{On the ring of a complex and the combinatory theory of integration}.
\item[2.] P.S.~ Aleksandrov (USSR), \textit{Some problems in the set-theoretic topology}.
\item[3.] G.~Birkhoff (USA), \textit{Continuous groups and linear spaces}.
\item[4.]  N.N.~Bogolyubov and N.M.~Krylov (USSR), \textit{Metric transitivity and invariant measure in dynamical systems of nonlinear mechanics}.
\item[5.] K.~Borsuk (Poland), \textit{On spheroidal spaces}.
\item[6.] N.K.~Brushlinskii (USSR), \textit{On continuous mappings of spherical manifolds}.
\item[7.] E. \v{C}ech (Czecho-Slovakia), \textit{Accessibility and homology}.
\item[8.] E. \v{C}ech (Czecho-Slovakia), \textit{Betti groups with different coefficient groups}. 
\item[9.] St.~Cohn-Vossen (USSR), \textit{Topological questions of differential geometry in the large}.
\item[10.] D. van Dantzig (Netherlands), \textit{Topological algebra}.
\item[11.] V.A.~Efremovich (USSR), \textit{On topological types of affine mappings}.
\item[12.] H.~Freudenthal (Netherlands), \textit{On topological approximations of spaces}.
\item[13.] I.I.~Gordon (USSR), \textit{On the intersection invariants of a complex and its residual space}.
\item[14.] P.~Heegaard (Norway), \textit{On the four-color problem}.
\item[15-16.] H.~Hopf (Switzerland), \textit{New 
research on $n$-dimensional manifolds}. Two~talks.
\item[17.] W.~Hurewicz (Netherlands), \textit{Homotopy and homology}.
\item[18.] E.R.~van Kampen (USA), \textit{On the structure of a compact group}.
\item[19.] A.N.~Kolmogorov (USSR), \textit{Homology rings in closed sets}.
\item[20.] K.~Kuratowski (Poland), \textit{On projective sets}.
\item[21.] S.~Lefschetz (USA), \textit{On locally connected sets}.
\item[22.] A.A.~Markov jr. (USSR), \textit{On the free equivalence of the closed braids}.
\item[23.] S.~Mazurkiewicz (Poland), \textit{On existence of non-decomposable continua in the sets of dimension $\geq 2$}.
\item[24.] V.V.~Nemytskii (USSR), \textit{Unstable dynamical systems}.
\item[25.] J.~von Neumann (USA), \textit{Integration theory in continuous
groups}.
\item[26-27.] J.~Nielsen (Denmark), Two talks on continuos surface mappings. 
\item[28.] G.~N{\"o}beling (Germany), \textit{On the triangulability of varieties and
main conjectures of combinatorial topology}.
\item[29.] L.S.~Pontryagin (USSR), \textit{Topological properties of compact
Lie groups}.
\item[30.] G.~de Rham (Switzerland), \textit{On new Reidemeister's topological invariants}.
\item[31.] G.~de Rham (Switzerland), \textit{Topological aspect of the theory of multiple integrals}.
\item[32.] J.A.~R\'{o}\.{z}a\'{n}ska (USSR), \textit{On continuous mappings of elements}.
\item[33.] J.~Schauder (Poland), \textit{Some applications of the topology of functional spaces}.
\item[34.] W.~Sierpi\'{n}ski (Poland), \textit{On contionuous mappings of sets}.
\item[35.] W.~Sierpi\'{n}ski (Poland), \textit{On transformations of sets by the Baire functions}.
\item[36.] W.~Sierpi\'{n}ski (Poland), \textit{On a projective set of the second class}.
\item[37.] P.A.~Smith (USA), \textit{Transformations of period two}.
\item[38.] M.H.~Stone (USA), \textit{Mappings theory in general topology}.
\item[39.] A.W.~Tucker (USA), \textit{On discrete spaces}.
\item[40.] A.N.~Tikhonov (USSR), \textit{On invariant points of continuous mappings of bicompact spaces}.
\item[41.] A.~Weil (France), \textit{Topological demonstration of the Cartan theorem}.
\item[42.] A.~Weil (France), \textit{The families of curves on the torus}.
\item[43.] H.~Whitney (USA), \textit{Topological properties of differentiable manifolds}.
\item[44.] H.~Whitney (USA), \textit{Sphere-spaces}.
\end{itemize}
\end{small}
\vspace{0.2cm}

Thus, we can see that the total number of talks actually  listed in \cite{LL57} equals 44. However, the authors of \cite{LL57} wrote about 45 talks including 13 made by mathematicians from USSR,\footnote{Notice that Stefan Cohn-Vossen is mentioned among speakers from Soviet Union. He was barred from lecturing in Cologne in 1933 under Nazi racial legislation. After a one-year period in Switzerland, he immigrated to USSR at the end of 1934. Up to the death from pneumonia in 1936, he worked at the Steklov Mathematical Institute, first in Leningrad, then in Moscow. See in this relation the speech of Cohn-Vossen at the beginning of his work in Leningrad \cite{CV}.} 
10 from USA, 7 from Poland, 4 from Switzerland, 3 from Netherlands, 2 from France, 2 from Czecho-Slovakia, 2 from Denmark, 1 from Germany, and 1 from Norway. We guess that Lapko and Lyusternik calculated N.N.~Bogolyubov and N.M.~Krylov as separate speakers despite of the fact that they made  a joint presentation.

The authors of the survey \cite{LL57} based heavily on the official Conference Report \cite{Al36}. Notice that  P.S.~Aleksandrov claims in \cite{Al36} about only 43 talks. Moreover, his report actually contains information on  42 presentations. Namely, the presentation of Gordon (see talk No 13 in the above-given list) and the second talk of Whitney (see No 44) are unmentioned in \cite{Al36}. Meanwhile, these both talks are presented in the official Proceedings published in \cite{MSb36}.\footnote{The talk of I.I.~Gordon is also mentioned in the report of A.W.~Tucker  \cite{Tuc35} and in the reminiscences of H.~Hopf \cite{H66}, while H.~Whitney in his historical note \cite{Whi88} points out that he gave two talks at the Moscow conference.}

Consider these Proceedings in more detail. It contains information about 41 talks. Indeed, the complete texts or extended abstracts are given for the talks Nos 2-5,  10, 12, 14, 19-21, 25, 30, 33-34, 37-38, and 42-44 (see the list of talks). The short abstracts are provided for presentations Nos 7, 17-18, 23, 32, 35-36, 39, and 41. For talks Nos 8, 13, 15-16, 22, 24, 26-29, and 40 there are references either to published papers or to ones in press, while for talks No 1 and No 9 only the titles are 
given.\footnote{Nevertheless, there is a reference to the ``abstract'' of talk No 1 in \cite{T17}.} There are no indications in \cite{MSb36} on presentations Nos 6, 11, and 31.

 Notice that the content of presentation No 1 is given in \cite{A35a} and \cite{A35b}; see also \cite{A36}. The contents of talks No 11 and No 31 are provided in \cite{E35} and \cite{dR36}, respectively. We guess  that talk No 9 is based on the paper \cite{CV36}. Unfortunately, we could not find any information about the content of talk No~6.
 
Observe also that the titles of some talks in the Proceedings \cite{MSb36} differ from those listed in the Conference Report \cite{Al36} and the survey~\cite{LL57}.\footnote{We emphasize that as the conference talks as items in Proceedings were presented in one of the four languages: English, French, German, and Russian. For the reader's convenience we give all the titles in English.} For instance, the papers of H.~Freundenthal, J.~von Neumann, and M.H.~Stone in \cite{MSb36} are titled \textit{Expansion of spaces and groups}, \textit{The uniqueness of Haar's measure}, and \textit{Applications of Boolean algebras to topology}, respectively. Also, in \cite{MSb36}, the exact title \textit{Topological invariants of the classes of surface mappings} of J.~Nielsen's presentation is provided.

Speaking about the scientific impact of the Conference, the authors of all reminiscences consider them outstanding. Many presentations contained the results that later became world-known named theorems. However, even within this framework, the birth of cohomology theory given in talks of Kolmogorov, Alexander and Gordon occupies a special place. 

Notice that Alexander and Kolmogorov have made (simultaniously and independently) the identical construction of the cohomology ring while Gordon's definition was slightly different.\footnote{Later Freundenthal proved isomorphism of these constructions.}

The second main import of the Conference, in Whitney's opinion \cite{Whi88}, was the connection between homotopy  and homology theories given in the talk of W.~Hurewicz. This became a basic result of algebraic topology. According to Whitney, Alexander, E.~\v{C}ech and D.~van Dantzig ``\dots also said that they had considered or actually used the definition of Hurewicz''.

The third remarkable result described in \cite{Whi88} in detail was an introduction of the Stiefel--Whitney classes presented in the talks of Hopf (who discussed the results of his PhD student E. Stiefel) and Whitney.\footnote{In \cite{Whi88}, Whitney explained that he decided to give ``two shorter talks'' under influence of Hopf's presentation.}

Some outcomes of the Conference have made a huge impact on neighboring fields of mathematics. Here we mention the celebrated Bogolyubov--Krylov theorem (the dynamical systems theory), the Leray--Schauder fixed point principle (nonlinear partial differential equations) and the Markov theorem (the braid theory).

Let us compare the factual list of talks with that in the scientific program of the Conference. This program was found by A.N.~Shiryaev  in the personal archive of A.N.~Kolmogorov during preparation of the three-volume Kolmogorov Commemorative Edition; now it is available at \cite{K03}. 

The tentative list of speakers in \cite[pp.590--593]{K03} does not contain the names of Efremovich and Sierpi\'{n}ski. Also, only one talk No 43 of Whitney was announced. On the other hand, the participation of A.W.~Iwanowski (USSR), Ch.H.~M{\"u}ntz (USSR), M.H.A.~Newman (UK), and K.~Reidemeister (Germany) was expected. Further, one more talk of Bogolyubov was planned, and the second talk of Hopf had another title.\footnote{We emphasize that the presumed number of talks in the program is also equal to 44.}

Besides the scientific presentations, the Conference program contains the special meeting of the Moscow Mathematical Society (September 5th) dedicated to the memory of the outstanding mathematician Emmy N\"{o}ther (1882--1935), whose works highly influenced the development of topology (see Aleksandrov's speech at this meeting \cite{Al36a}).

It should be noted that in September 1935, \textit{Bulletin of the American Mathematical Society} published an announcement of the Conference with the names of American mathematicians who expected to attend and take an active part (Section \textit{Notes}, see \cite[p.615]{Bul35}). Besides 9 persons already mentioned above, this list contains the names of D.V.~Widder and O.~Zariski. 

The presence of Widder at Moscow conference is confirmed by his own reminiscences \cite[p.82]{Wid88}, while G.~Birkhoff mentioned in the memoir \cite[p.45]{Bi89} that Widder attended the Conference ``informally'' (without any talk). 

As for Zariski, this is a more complicated case. In \cite{Tuc35} Tucker claimed that Zariski gave a conference talk on algebraic geometry.  This apparent contradiction can be resolved by using Zariski's biography \cite{P09}. We can read at \cite[p.60]{P09} that ``...Weil and Zariski (who had not been an official delegate) were invited by the mathematical department at the University of Moscow to give a series of lectures...''

In \cite[p.331]{Al80}, P.S.~Aleksandrov mentioned  B.~Knaster (Poland) among the speakers. This definitely contradicts  all the other sources. Moreover, we did not find any documentary evidence of Knaster's presence at the Conference.\footnote{Prof. J. Mioduszewski, a former student of B. Knaster, guesses that Knaster could attend the 
Conference privately.}

We also mention the report of K.~Borsuk \cite{B36} where the number of speakers  by countries (except USSR) was pointed out. These data correspond to ours given above.

Surely, there were many conference attendees without any talk including the well-known mathematicians. For instance, present at the Opening Ceremony was Academician N.N.~Luzin  (see, e.g., the report \cite{Z35} of K.~Zarankiewicz, who was a member of the  official Polish delegation). 

L.V.~Kantorovich (later Nobel Prize winner) says in his reminiscences \cite[p.250]{Ka87}:  ``... I even thought about including a supplementary report at the forthcoming conference on topology in Moscow in September, but apparently <...>  became quite seriously ill. Although I attended the conference itself, I was not strong enough to prepare a contribution at an appropriate level.''

\begin{figure*}
\centering
  \includegraphics[width=0.75\textwidth]{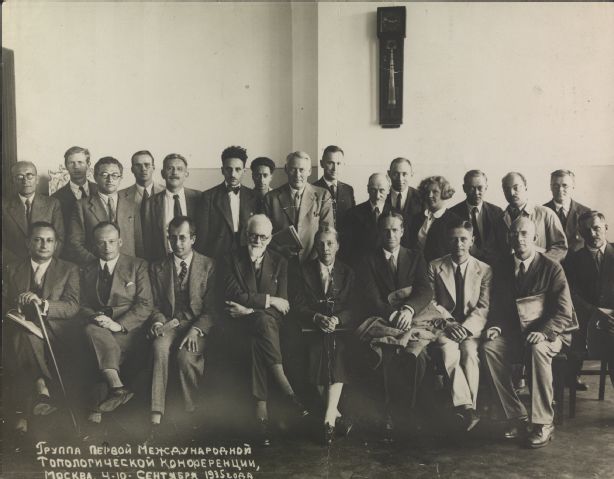}
\caption{First International Topology Conference in Moscow}
\label{fig:group}       
\end{figure*}

Now we would like to present the biggest curiosity related to the Conference. In \cite[pp.96--97]{Whi88}, Whitney  provided a  group photo of the participants  with (incomplete) description (see Fig.~\ref{fig:group}).

In particular, a man, the 5th from the right in the last row on Fig.~\ref{fig:group},  is called J.D.~Tamarkin.  This description was repeated in Whitney's Collected Papers \cite{Whi92} and in some other publications, e.g., in the monograph \cite{J99}  and in the online archive of ETH Z{\"u}rich \cite{ETH}. Moreover, the image copied from this photo is given as a portrait of Tamarkin in Wikipedia.\footnote{except the Russian version.}

We claim that this is a mistake. Indeed, Tamarkin was a prominent analyst, but he never worked in topology and related fields (see \cite{Hi47}). 

Next, Tamarkin left  the Soviet Union in December 1924 illegally,\footnote{As a former member the {\it Menshevik} (social democratic) party Tamarkin was a ``fair game'' for the Soviet secret police.} crossing the border through the frozen  Chudskoe Lake with smugglers (see \cite{Sha05}). Thus, it is highly unlikely that the Soviet authorities  allowed the emigrant Tamarkin to enter the country. On the other hand, Tamarkin, who feared arrest as far back as 1924, would not have travelled to  the USSR in 1935. 

\begin{figure}[htbp]
\begin{minipage}[b]{0.5\textwidth}
\centering
\includegraphics[width=5.6cm]{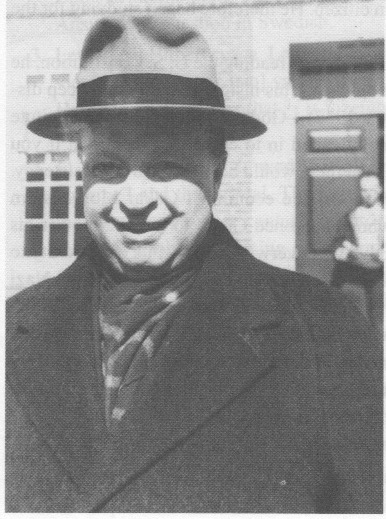}
\caption{J.D.~Tamarkin}
\label{Fig2}
\end{minipage}
\hfill
\begin{minipage}[b]{0.47\textwidth}
\centering
\includegraphics[width=4.7cm]{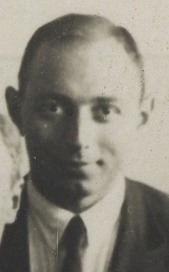}
\caption{}
\label{Fig3}
\end{minipage}
\end{figure}

\begin{figure}[h!]
\centering
\includegraphics[width=0.75\textwidth]{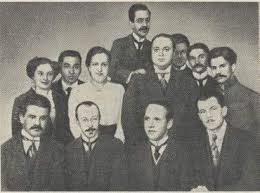}
\caption{St. Petersburg. J.D.~Tamarkin  in the middle of the 2nd row }
\label{seminar}
\end{figure}

\begin{figure}[h!]
\centering
\includegraphics[width=0.6\textwidth]{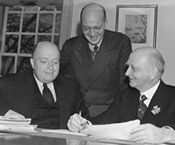}
\caption{Brown University. J.D.~Tamarkin 1st from the left }
\label{Brown}
\end{figure}

Notice also that Tamarkin was 47 years old in 1935 and he was  corpulent (see \cite[p.115]{Tol96} and Fig.~\ref{Fig2}), whereas the man in the group photo (Fig.~\ref{Fig3}) looks much younger and thinner. In addition, if one looks at two well-known photos of Tamarkin (Fig.~\ref{seminar} of 1914  and Fig.~\ref{Brown} of 1941) then the mistake in \cite{Whi88} becomes obvious.

Finally, as we have found out, Tamarkin  participated in the meeting of the American Mathematical Society at Ann Arbor (September 10--13, 1935) and, in particular, made a speech at a joint dinner of the AMS, the Mathematical Association of America, and the Institute of Mathematical Statistics on September 12 (see \cite[p.755]{Ri35} and \cite[p.32]{MAA}). Obviously, at that time it was impossible to get from Moscow to Ann Arbor in several days.

To clarify who is this person on the group photo, we need to turn to the history of the  Moscow school of topology. The leader of Moscow topologists in the thirties was P.S.~Aleksandrov. One of his favorite and talented students was Lev A. Tumarkin (1904--1974), who established some fundamental results in dimension theory.
In 1935 he became the Dean of the  Mechanics and Mathematics Faculty at Moscow State University. Thus, it seems rather natural to see Tumarkin among the conference participants as a hospitable host, even without a presentation. 

\begin{figure}[h!]
\begin{minipage}[b]{0.49\textwidth}
\centering
\includegraphics[width=5.2cm]{Non-Tamarkin.jpg}
\caption{}
\label{Fig4}
\end{minipage}
\hfill
\begin{minipage}[b]{0.49\textwidth}
\centering
\includegraphics[width=5.6cm]{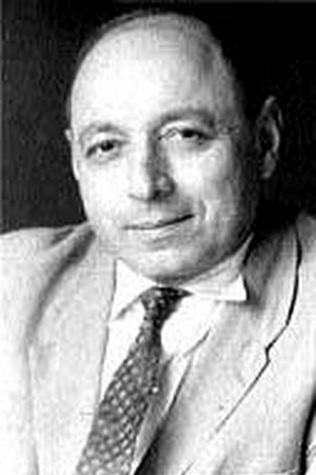}
\caption{L.A.~Tumarkin}
\label{Fig5}
\end{minipage}
\end{figure}

Unfortunately, we do not have a Tumarkin's photo of that time. But comparing a man from the group photo with  the portrait of Tumarkin at the age 60 (see Fig.~\ref{Fig5}), we see that it could very well be the same person. Therefore, we conclude that H.~Whitney was misled by the similarity of surnames T$\!$\textbf{\textit{a}}markin and 
T$\!$\textbf{\textit{u}}markin.

In conclusion, we notice that in Whitney's description of the Conference photo  there were also three unknowns. Prof. L.~Maligranda \cite{M18} identified these people as geometers from Kharkov: D.M.~Sintsov, M.A.~Nikolaenko, and P.A.~Solov'yov.\footnote{Some information about these persons can be found in \cite{Kha}.} Now, we can provide  a complete description. \medskip

The \ first \ row \ (sitting, \ left \ to \ right): \ K.~Kuratowski (1896--1980), \ J.P.~Schauder (1899--1943), \ S.~Cohn-Vossen (1902--1936), \ P.~Heegaard (1871--1948), \ J.A.~R\'{o}\.{z}a\'{n}ska\footnote{Julia A. R\'{o}\.{z}a\'{n}ska was a topologist, a disciple of P.S. Aleksandrov, Associate Professor at Moscow University, Moscow Mathematical Society member.} (1901--1967), J.W.~Alexander (1888--1971), H.~Hopf (1894--1971), P.S.~Aleksandrov (1896--1982), P.A.~Solov'yov (1906--1993).

The second row (standing, left to right): \ \ E.~\v{C}ech (1893--1960),\ \  H.~Whitney (1907--1989),\  \ \ K.~Zarankiewicz (1902--1959), A.W.~Tucker (1905--1995), \ S.~Lefschetz (1884--1972), \ H.~Freudenthal (1905--1990), \ F.I.~Frankl (1905--1961), J.~Nielsen (1890--1959), \ \ K.~Borsuk (1905--1982), \ \ D.M.~Sintsov (1867--1946), \ L.A.~Tumar\-kin (1904--1974), \  \ 
 M.A.~Nikolaenko\footnote{Maria A. Nikolaenko was a geometer, Associate Professor at Kharkov University.} (1905--1988),\ \  V.V.~Stepanov (1889--1950), E.R.~van Kampen (1908--1942), A.N.~Tikhonov (1906--1993).

\subsection*{Acknowledgements}
\begin{small} We are grateful to Lech Maligranda who pointed out earlier unknown names from the Conference photo and provided us with the paper \cite{Z35}. We also would like to thank Jerzy Mioduszewski, Roman Duda and Galina Smirnova for useful discussions, and Sem\"en S.~Kutateladze and Ignat Soroko for valuable comments. The first author was partly supported by the ``RUDN University Program 5-100''.
\end{small}

\Addresses


\begin{thebibliography}{MSb36}
%
%

\bibitem{A35a}
J.W. Alexander.
\newblock On the chains of a complex and their duals.
\newblock {\em Proc. Natl. Acad. Sci. USA} 21 (1935), 509--511.

\bibitem{A35b}
J.W. Alexander.
\newblock On the ring of a compact metric space.
\newblock {\em Proc. Natl. Acad. Sci. USA} 21 (1935), 511--512.

\bibitem{A36}
J.W. Alexander.
\newblock On the connectivity ring of an abstract space.
\newblock {\em Ann. of Math. (2)} 37(3) (1936), 698--708.

\bibitem{Al36}
P.S. Aleksandrov.
\newblock First international congress on topology.
\newblock {\em Uspekhi Mat. Nauk} 1 (1936), 260--262.
\newblock [Russian].

\bibitem{Al36a}
P.S. Aleksandrov.
\newblock In memory of Emmy Noether, 
\newblock {\em Uspekhi Mat. Nauk} 2 (1936), 255--265.
\newblock [Russian]. English translations in: A. Dick.
Emmy Noether. 1882--1935, Birkh\"{a}user, 1981, pp. 153--179.

\bibitem{Al80}
P.S. Aleksandrov.
\newblock Pages from an autobiography. {II}.
\newblock {\em Russian Math. Surveys} 35(3) (1980), 315--358.

\bibitem{Bi89}
G.~Birkhoff.
\newblock Mathematics at {H}arvard, 1836--1944.
\newblock In {\em A century of mathematics in {A}merica, {P}art {II}}, volume~2
  of {\em Hist. Math.}, pp. 3--58. Amer. Math. Soc., Providence, RI, 1989.

\bibitem{B36}
K.~Borsuk.
\newblock Mi{\k{e}}dzynarodowa {K}onferencja {T}opologiczna w {M}oskwie.
\newblock {\em Wiadomo{\'s}ci Matematyczne}, 41 (1936), 134--137.
\newblock [Polish].

\bibitem{Bul35}
Notes.
\newblock {\em Bull. Amer. Math. Soc.}, 41(09) (1935), 615.

\bibitem{CV}
	St. Cohn-{V}ossen.
\newblock A speech of {S}. {C}ohn-{V}ossen at the beginning of {S}. {C}ohn-{V}ossen's {U}{S}{S}{R} career. 
 http://www.mi.uni-koeln.de/ Geschichte/COHN-VOSSEN,Stefan/PersonalRemarks.pdf

\bibitem{CV36}
St. Cohn-{V}ossen.
\newblock Bending of surfaces in the large.
\newblock {\em Usp. Mat. Nauk} 1 (1936), 33--76.
\newblock [Russian].

\bibitem{dR36}
G.~de~Rham.
\newblock {Relations entre la topologie et la th\'eorie des int\'egrales
  multiples.}
\newblock {\em Enseign. Math.} 35 (1936), 213--228. 
\newblock [French]. 

\bibitem{E35}
V.A. Efremovich.
\newblock Topological classification of affine mappings of the plane.
\newblock {\em Sb. Math.} 42(1) (1935), 23--29 [Russian], 30--36 [German].

\bibitem{ETH}
E{T}{H}-{B}ibliothek {Z}{\"u}rich, {B}ildarchiv, {P}ortr{\_}07561, 
  http://doi.org/ 10.3932/ethz-a-000491996

\bibitem{Hi47}
E.~Hille.
\newblock Jacob {D}avid {T}amarkin---his life and work.
\newblock {\em Bull. Amer. Math. Soc.} 53 (1947), 440--457.

\bibitem{H66}
H.~Hopf.
\newblock Einige pers\"onliche {E}rinnerungen aus der {V}orgeschichte der
  heutigen {T}opologie.
\newblock Centre Belge Rech. Math., Colloque Topologie, Bruxelles 1964 
  (1966), 440--457.
\newblock [German].

\bibitem{J99}
I.M. James, editor.
\newblock {\em History of topology}.
\newblock North-Holland, Amsterdam, 1999.

\bibitem{Ka87}
L.V. Kantorovich.
\newblock My journey in science (proposed report to the {M}oscow {M}athematical
  {S}ociety).
\newblock {\em Russian Math. Surveys} 42(2) (1987), 233--270.


\bibitem{Kha}
On the {K}harkov geometrical school. 
http://www.univer.omsk.su/ omsk/Sci/HkGS/hkgs3.html
\newblock [Russian].


 
\bibitem{LL57}
A.F. Lapko and L.A. Lyusternik.
\newblock Mathematical sessions and conferences in the {USSR}.
\newblock {\em Uspehi Mat. Nauk (N.S.)}, 12(6(78)) (1957), 47--130.
\newblock [Russian].

\bibitem{M18}
L.~Maligranda.
\newblock Personal communication, 2018.

\bibitem{MSb36}
The {F}irst {I}nternational {T}opological {C}onference, {S}eptember 4--10, 1935.
\newblock {\em Sb. Math.} 1 (43)(5), 1936.

\bibitem{P09}
C.~Parikh.
\newblock {\em The unreal life of {O}scar {Z}ariski}.
\newblock Springer, New York, 2009.
\newblock With a foreword by David Mumford, Reprint of the 1991 original.

\bibitem{Ri35}
R.G.D. Richardson.
\newblock The summer meeting at {A}nn {A}rbor.
\newblock {\em Bull. Amer. Math. Soc.}, 41(11) (1935), 753--763.

\bibitem{MAA}
K.~Ross and J.~Tattersall.
\newblock Meetings of the {M}{A}{A}.
  https://www.maa.org/about-maa/maa-history/celebrating-the-centennial/maa-meetings

\bibitem{Sha05}
T.~Shaposhnikova.
\newblock Three high-stakes math exams.
\newblock {\em Math. Intelligencer} 27(3) (2005), 44--46.

\bibitem{K03}
A.N. Shiryaev, editor.
\newblock {\em Kolmogorov. {C}ommemorative edition in three books. {B}ook 2.
  {T}he braid of these running lines ... {S}elected passages from the
  correspondence between {A}.{N}. {K}olmogorov and {P}.{S}. {A}leksandrov}.
\newblock Moskva: Fizmatlit ``Nauka'', 2003.
\newblock [Russian].

\bibitem{S17}
G.S. Smirnova.
\newblock The {F}irst {I}nternational {T}opology {C}onference. {M}oscow 1935.
\newblock In {\em Proceedings of the International Conference Analytical and
  Computational Methods in Probability Theory and its Applications}, pp.
  187--191. RUDN, Moscow, 2017.
\newblock [Russian].

\bibitem{T17}
F.~Toenniessen.
\newblock {\em Topologie. Ein Lesebuch von den elementaren Grundlagen bis zur
  Homologie und Kohomologie}.
\newblock Heidelberg: Springer Spektrum, 2017.
\newblock [German].

\bibitem{Tol96}
E.~Tolsted.
\newblock Reminiscences about {P}rofessor {Y}a.{D}. {T}amarkin.
\newblock {\em Istor.-Mat. Issled. (2)} (1(36), part 2) (1996), 108--118.
\newblock Publication translation and notes by N.S. Ermolaeva. 
\newblock [Russian].

\bibitem{Tuc35}
A.W.~Tucker.
\newblock The topological congress in {M}oscow.
\newblock {\em Bull. Amer. Math. Soc.} 41(11) (1935), 764.

\bibitem{Whi88}
H.~Whitney.
\newblock Moscow 1935: topology moving toward {A}merica.
\newblock In {\em A century of mathematics in {A}merica, {P}art {I}}, volume~1
  of {\em Hist. Math.}, pp. 97--117. Amer. Math. Soc., Providence, RI, 1988.

\bibitem{Whi92}
H.~Whitney.
\newblock {\em The collected papers of {H}assler {W}hitney}, volume~I of {\em
  Contemp. Mathematicians}.
\newblock Boston, MA etc.: Birkh\"auser, 1992.

\bibitem{Wid88}
D.V. Widder.
\newblock Some mathematical reminiscences.
\newblock In {\em A century of mathematics in {A}merica, {P}art {I}}, volume~1
  of {\em Hist. Math.}, pp. 79--83. Amer. Math. Soc., Providence, RI, 1988.

\bibitem{Z35}
K.~Zarankiewicz.
\newblock Mi{\k{e}}dzynarodowa {K}onferencja {T}opologiczna w {M}oskwie.
\newblock {\em Mathesis Polska} 10(5--6) (1935), 114--116.
\newblock [Polish].
\end{thebibliography}
\end{document}